\documentstyle{amlts}
\begin{document}
\annalsline{158}{2003}
\received{May 8, 2002}
\startingpage{1067}
\def\bye{\end{document}}
 \font\tenrm=cmr10
\def\ritem#1{\item[{\rm #1}]}
\catcode`\@=11
\font\twelvemsb=msbm10 scaled 1100
\font\tenmsb=msbm10
\font\ninemsb=msbm10 scaled 800
\newfam\msbfam
\textfont\msbfam=\twelvemsb  \scriptfont\msbfam=\ninemsb
  \scriptscriptfont\msbfam=\ninemsb
\def\msb@{\hexnumber@\msbfam}
\def\Bbb{\relax\ifmmode\let\next\Bbb@\else
 \def\next{\errmessage{Use \string\Bbb\space only in math
mode}}\fi\next}
\def\Bbb@#1{{\Bbb@@{#1}}}
\def\Bbb@@#1{\fam\msbfam#1}
\catcode`\@=12

 \catcode`\@=11
\font\twelveeuf=eufm10 scaled 1100
\font\teneuf=eufm10
\font\nineeuf=eufm7 scaled 1100
\newfam\euffam
\textfont\euffam=\twelveeuf  \scriptfont\euffam=\teneuf
  \scriptscriptfont\euffam=\nineeuf
\def\euf@{\hexnumber@\euffam}
\def\frak{\relax\ifmmode\let\next\frak@\else
 \def\next{\errmessage{Use \string\frak\space only in math
mode}}\fi\next}
\def\frak@#1{{\frak@@{#1}}}
\def\frak@@#1{\fam\euffam#1}
\catcode`\@=12

\newcommand{\PSI}{\vert\psi\rangle}
\newcommand{\PHI}{\vert\phi\rangle}
\newcommand{\D}{{\cal D}}
\newcommand{\pd}{\partial}
\newcommand{\eps}{\varepsilon}
\newcommand{\suli}{\sum\limits}
\newcommand{\inli}{\int\limits}
\newcommand{\bh}{{b/2}}
\newcommand{\aob}{\left(\frac{a}{b}\right)}
\newcommand{\V}{V}
\newcommand{\al}{\alpha}
\newcommand{\xij}{|x_i-x_j|}
\newcommand{\half}{\mbox{$\frac{1}{2}$}}
\newcommand{\E}{\cal {E}}
\newcommand{\FF}{\cal {F}}
\newcommand{\F}{F}
\newcommand{\G}{G}
\newcommand{\rmax}{\rho_{\al,{\rm max}}}
\newcommand{\rmin}{\rho_{\al,{\rm min}}}
\newcommand{\R}{{\Bbb R}}
\newcommand{\real}{{\rm Re}}
\newcommand{\imag}{{\rm Im}}
\newcommand{\e}{\widetilde{e}}
\newcommand{\vv}{\widetilde{v}}
\newcommand{\RR}{\widetilde{R}}
\newcommand{\as}{\widetilde{a}}
\newcommand{\beq}{\begin{equation}}
\newcommand{\eeq}{\end{equation}}
\newcommand{\const}{{\rm const. \,}}
\newcommand{\supp}{{\rm supp \,}}

\title{Poincar\'e  inequalities\\ in punctured domains} 
\shorttitle{Poincar\'e  inequalities in punctured domains} 

  \acknowledgements{The work of the first named author was  partially
supported by U.S. National Science Foundation
grant PHY 98-20650.  The   second named author is an Erwin Schr\"odinger Fellow,
supported by the Austrian Science Fund.\hfill\break
\hglue19pt  \raise 1pt\hbox{\tiny\copyright}2003 by the authors.  This paper may be reproduced, in its entirety, for
noncommercial purposes.}
 \twoauthors{Elliott H.~Lieb, Robert Seiringer,}{Jakob
Yngvason}
\institutions{Departments of Mathematics and Physics, Princeton University, Princeton,
NJ\\
{\eightpoint {\it E-mail address\/}: lieb@math.princeton.edu}
\\
\vglue6pt
Department of Physics, Princeton University, Princeton, NJ\\
{\eightpoint {\it E-mail address\/}:  rseiring@princeton.edu}
\\
\vglue6pt
Institute for Theoretical Physics, University of Vienna,  Vienna, Austria
\\
{\eightpoint {\it E-mail address\/}: yngvason@thor.thp.univie.ac.at}}

\centerline{\bf Abstract}
\vglue6pt

The classic Poincar\'e inequality bounds the $L^q$-norm of a function
$f$ in a bounded domain $\Omega
\subset \R^n$ in terms of some $L^p$-norm of its gradient in $\Omega$.
We generalize this in
two ways: In the first generalization we remove a set $\Gamma$ from
$\Omega$ and concentrate our attention on $\Lambda =
\Omega \setminus \Gamma$. This new domain might not even be connected and
hence no Poincar\'e inequality can generally hold for it, or if it
does hold it might have a very bad constant.  This is so even if the
volume of $\Gamma$ is arbitrarily small.  A Poincar\'e inequality {\it
does hold}, however, if one makes the additional assumption that $f$
has a finite $L^p$ gradient norm on the whole of $\Omega$, not just on
$\Lambda$. The important point is that the Poincar\'e inequality thus
obtained bounds the $L^q$-norm of $f$ in terms of the $L^p$ gradient
norm on $\Lambda$ (not $\Omega$) plus an additional term that goes to
zero as the volume of $\Gamma$ goes to zero.  This error term depends
on $\Gamma$ only through its volume. Apart from this additive error
term, the constant in the inequality remains that of the `nice' domain
$\Omega$.  In the second generalization we are given a vector field
$A$ and replace $\nabla $ by $\nabla +i A(x)$ (geometrically, a
connection on a $U(1)$ bundle). Unlike the $A=0$ case, the infimum of
$\Vert (\nabla +i A)f \Vert_p $ over all $f$ with a given $\Vert f
\Vert_q$ is in general not zero. This permits an improvement of the
inequality by the addition of a term whose sharp value we derive. We
describe some open problems that arise from these
generalizations.

\vglue-4pt
\section{Introduction}
\vglue-4pt

The simplest Poincar\'e inequality refers to a bounded, connected domain
$\Omega \subset \R^n$, and  a function $f\in
L^2(\Omega)$ whose distributional gradient is also in $L^2(\Omega)$
(namely, $f\in
W^{1,2}(\Omega)$). While it is false that there is a finite constant~$S$,
depending only on $\Omega$, such that
\begin{equation}\label{l2}
\int_\Omega |f|^2 \leq  S \int_\Omega |\nabla f|^2
\end{equation}
for all $f$, such an inequality {\it does}
hold if we impose the additional condition that $\int_\Omega f =0$. The
constant $S$ depends on $\Omega$, but it is independent of $f$. In fact,
$1/S$ is the second eigenvalue of the Laplacian in $\Omega$ with Neumann
boundary conditions. This is merely a consequence of Bessel's inequality.

Some simple generalizations of (\ref{l2}) are well known. One involves replacing the
condition $\int_\Omega f =0$ by the condition $\int_\Omega fg =0$,
where $g$ is any $L^2(\Omega)$ function that is not orthogonal to the
lowest Neumann eigenfunction of the Laplacian, i.e., $\int_\Omega g \neq 0$.
Another involves  replacing $L^2$ by $L^p$
for $1\leq p\leq \infty$ on both sides of (\ref{l2}) . Finally,
by Sobolev's inequality, the $L^p$-norm of $f$ can be replaced
by a suitable $L^q$-norm with $q >p$. The situation is summarized in
(see, e.g., \cite[Thms.~8.11 and~8.12]{anal}) the following statement:

\proclaimtitle{Standard Poincar\'e inequalities for
$W^{1,p}(\Omega)$} \specialnumber{1}\proclaim{Theorem} 
\label{T1}
Let $\Omega \subset \R^n$ be a bounded{\rm ,} connected{\rm ,} open set with
the cone property. Let $1\leq q \leq \infty$ and let
$\max\{1,qn/(n+q)\} \leq p \leq
\infty$  if $q< \infty$ and $n<p\leq \infty$ if $q=\infty$.  Let $g$ be a
function in $L^{p'}(\Omega)${\rm ,} $p^\prime = p/(p-1)${\rm ,} such that $\int_\Omega
g =1 $.  Then there is a constant $S_{p,q}>0${\rm ,} which depends on
$\Omega,\, g,\, p, \, q${\rm ,} such that for any $f\in W^{1,p}(\Omega),$
\begin{equation} \label{poin}
\Big\Vert f- \int_{\Omega}fg \Big\Vert_{L^q(\Omega)} \leq S_{p,q} \Vert
\nabla f \Vert_{L^p(\Omega)} \  .
\end{equation}
\endproclaim

{\it Remarks.} The case  $q=np/(n-p) $ requires
the Sobolev inequality explicitly for the proof, and thus the
inequality can be called the Poincar\'e-Sobolev inequality in this
case.  The domain $\Omega$ is required to have the ``cone property"
(see, e.g., \cite{mazja}); i.e., each point of $\Omega$ is the vertex
of a spherical cone with fixed height and angle,   which is situated
in $\Omega$. Note that the constants $S_{p,q}$ depend not only on the
volume of $\Omega$, but also on its shape.  The Poincar\'e inequality
is usually presented as (\ref{poin}) with $q=p$ and with
$g=1/|\Omega|$, where, in general, $|\cdot|$ denotes the Lebesgue
measure of a set in $\R^n$.  A generalization to $W^{m,p}(\Omega)$ for
$m>1$ is possible (see, e.g.,
\cite{anal}).
\vglue12pt

Now we turn to the two generalizations that concern us in this
paper. They were motivated by a treatment of the quantum mechanical
many-body problem, specifically, the proof of Bose-Einstein
condensation in a physically realistic model \cite{LS}. Further
developments required a version of the Poincar\'e inequality in which
$\nabla$ is replaced by a connection on a $U(1)$ bundle, namely
$\nabla \to \nabla + iA$, \pagebreak where $A$ is a vector
field. The generalization to this situation leads to a proof of
superfluidity   for the same quantum-mechanical system
\cite{lsy02}. Note that for $n=1$ the bundle is trivial, since the 
vector field $A$ can be eliminated by a unitary transformation, namely
$f(x)\to f(x)\exp(-i\int_z^x A(y)dy)$ for some $z\in\Omega$.

\vglue12pt
Our main result is Theorem  \ref{T4}, which contains the two generalizations,
and which we now describe in detail.

The first generalization concerns the following obstruction to the use of the
Poincar\'e inequality (\ref{poin}): Let us remove a small set $\Gamma$ from
$\Omega$ and concentrate our attention on $\Lambda = \Omega \setminus
\Gamma$. This new domain might not even be connected and hence no
Poincar\'e inequality can generally hold for $\Lambda$, no matter how
small $|\Gamma|$ might be. Even if $\Lambda $ is connected, the
constant $S_{p,q}$ could be very large, or even infinite.

A trivial example is to let $\Omega$ be a unit square in $\R^2$ and to
let $\Gamma$ be a thin annulus in $\Omega$ of outer radius 1/2 and
inner radius $1/2 - \varepsilon$. Take $g(x) = 1$.  We can take $f=1$
inside the disk of radius $1/2 - \varepsilon$ and $f=0$ elsewhere.
Thus, regardless of how small $\varepsilon$ may be, the right side of
(\ref{poin}), with $\Omega$ replaced by $\Lambda$, will be zero while
the left side is positive.  Another, perhaps more interesting example
is one in which $\Lambda$ is connected but fails to satisfy
Theorem~\ref{T1} because the cone property is absent. This can be
accomplished with a small $\Gamma$ that is topologically a ball, but
which has a sufficiently rough surface (see, e.g.,\break
\cite[\S 8.7]{anal}).

The smallness of $|\Gamma|$ cannot restore the Poincar\'e inequality
in $\Lambda$.  A generalized Poincar\'e inequality {\it does hold},
however, if one makes the additional assumption that $f$ has
an extension to a function with a finite
$L^p$ gradient norm on the whole of $\Omega$, not just on
$\Lambda$. The important point is that the Poincar\'e inequality thus
obtained bounds the $L^q$-norm of $f$ in terms of the $L^p$ gradient
norm on $\Lambda$ (not $\Omega$) plus an additional error
term. Furthermore -- and this is also important -- the effective
Poincar\'e inequality that holds for $\Lambda$ approaches that given
in (\ref{poin}) as the volume of $\Gamma$ tends to zero --- in a
manner that depends only on the fixed $L^p$ gradient
norm on $\Omega$ and on $|\Gamma|$, but not on its shape. There is
no regularity assumption on $\Gamma$.

In the second generalization $\nabla$ is replaced by $\nabla +i A(x)$
on the right side of (\ref{poin}), where $A:\Omega\to\R^n$ is some
given vector field. (For simplicity we assume that $A$ is a bounded,
measurable function, but a weaker condition will certainly suffice
(see, e.g., \cite[\S 7.20]{anal}).) 
We observe that hidden in (\ref{poin}) is a (nonlinear)
`eigenvalue' that happens to be zero, namely the smallest value of
$\Vert \nabla f\Vert _p$ given the value of $\|f\|_q$. Thus,
(\ref{poin}) states that if the right side of (\ref{poin}) is small
then $f$ must be close to the `lowest' eigenfunction, which happens to
be the constant function.  Our goal is to do something similar when
$A(x) \neq 0$, i.e.,  to show that if $\Vert (\nabla
+iA(x))f\Vert _p$ is close to the lowest $L^p\to L^q$ `eigenvalue'
defined in (\ref{defea2}) then $f$ must be close to the corresponding
`eigenfunction'.

For $1\leq q \leq \infty$ and
$\max\{1,qn/(n+q)\} < p \leq
\infty$ the `energy'
$E_A^{p,q}$ is defined by
\begin{equation}\label{defea2}
E_A^{p,q}= \inf \left\{ \frac{\Vert (\nabla +i A)f\Vert_{L^p(\Omega)}}{\Vert f\Vert_{L^q(\Omega)}} \, : \, f\in W^{1,p}(\Omega),\, f\neq 0 \right\} \ .
\end{equation}
The {\it ground state manifold} ${\cal 
M}_A^{p,q}$ is given by the set of minimizers of (\ref{defea2}), which
is nonempty as a consequence of Theorem~\ref{Tmagn} below.  Note
that, in general, this will not be a linear space. Also its dimension
can be greater than one, as we show in the appendix for the case
$p=q=2$.

When $A=0$ the ground state manifold ${\cal  M}_A^{p,q}$ is one-dimensional, spanned by the constant function. In this case the replacement of
$f$ by $f-\int fg$ as in (\ref{poin}) has the same qualitative effect as
restricting the inequality to functions $f$ whose $L^q(\Omega)$
distance to the constant function is bounded below by a fixed multiple of the
$L^q(\Omega)$ norm of $f$. The fixed multiple depends on $g$, of
course, but so does the constant $S_{p,q}$ appearing in (\ref{poin}).

For $A\neq 0$ the dimension of ${\cal  M}_A^{p,q}$ can be greater than $1$, and we adopt
the second viewpoint in this case. We obtain an inequality for
functions whose distance $d_A^q$ to ${\cal  M}_A^{p,q}$ exceeds a
certain value $\delta>0$, i.e.,
\begin{equation}\label{distanceq}
d_A^q(f):=\inf_{\phi\in {\cal  M}_A^{p,q}} \|f-\phi\|_{L^q(\Omega)} \geq
\delta \,
\|f\|_{L^q(\Omega)} \ .
\end{equation}

Before giving our main Theorem \ref{T4} on punctured domains with
$A$ \break fields it may be useful to state
the following theorem, which generalizes\break Theorem~\ref{T1} to the case
of $A$ fields alone. We consider here only the case
$p>\break \max\{1,qn/(n+q)\}$, and leave the case of the `critical'
$p=qn/(n+q)\geq 1$ as well as the case $p=1$ as open problems.

\proclaimtitle{Poincar\'e inequalities with vector fields} \specialnumber{2}\proclaim{Theorem}  
\label{Tmagn}
Let $\Omega \subset \R^n$ be a bounded{\rm ,} connected{\rm ,} open set with
the cone property. Let $1\leq q \leq \infty${\rm ,}
$\max\{1,qn/(n+q)\} < p \leq\infty${\rm ,}
 and let $0<\delta\leq 1$. Then there is a constant $S^{p,q}_\delta>0${\rm ,} which depends on
$\Omega,\,\delta ,\, p, \, q${\rm ,} such that for any $f\in
W^{1,p}(\Omega)$ with\break\vglue-11pt \noindent $d^q_A(f)\geq \delta \|f\|_{L^q(\Omega)},$
\begin{equation}\label{ineqpq}
\|f \|_{L^q(\Omega)} \leq
S^{p,q}_\delta \left[ \|(\nabla+iA) f\|_{L^p(\Omega)} - E_A^{p,q}
\|f\|_{L^q(\Omega)}
\right] \ .
\end{equation}
\endproclaim

Note that Theorem~\ref{Tmagn} implies in particular that ${\cal  M}^{p,q}_A$ is not the empty set.

\demo{{R}emark}
If ${\cal  M}_A^{p,q}$ is one-dimensional, spanned by the minimizer
$\phi_A$, which we assume to be normalized by
$\|\phi_A\|_{L^q(\Omega)}=1$, we can go back to our original
formulation and take $g$ to be an $L^{q'}(\Omega)$ function satisfying
$\int_{\Omega} g\phi_A =1$ (compare with Theorem~\ref{T1}). The
corresponding generalization of (\ref{poin}) is then
\begin{equation} \label{l5}
\Big\| f(\cdot) -\phi_A(\cdot) \mbox{$\left[\int_\Omega fg \right]$}
\Big\|_{L^q(\Omega)}
 \leq \widetilde
S^{p,g}_{g} \left[ \|(\nabla+iA) f\|_{L^p(\Omega)} - E_A^{p,q}
\|f\|_{L^q(\Omega)}
\right] \  ,
\end{equation}
which now holds for all $f\in W^{1,p}(\Omega)$. Here $\widetilde
S^{p,q}_g$ is some constant depending on $g$ (besides $p$, $q$ and
$\Omega$). For $q<\infty$ a possible choice for $g$ would be $g(x)=\overline
{\phi_A(x)} |\phi_A(x)|^{q-2}$ if $\phi_A(x)\neq 0$, and $g(x)=0$
otherwise.
\enddemo

The generalization to punctured domains is the following, which is our
main theorem.

\proclaimtitle{Poincar\'e inequalities in punctured domains}
\specialnumber{3} \proclaim{Theorem}  \label{T4}
Let $1\leq q\break\leq \infty $ and $\max\{1, qn/(n+q)\}< p\leq\infty$.  Let
$\Omega \subset \R^n$ be a bounded{\rm ,} connected{\rm ,} open set with the cone
property{\rm ,} and let $E_A^{p,q}$ and ${\cal  M}_A^{p,q}$ be as
explained above. Let $\Lambda\subset \Omega$ be a measurable subset of
$\Omega${\rm ,} $\Gamma=\Omega\setminus\Lambda${\rm ,} and let $0<\delta\leq
1$. For any $\eps>0$ there exists a positive constant $C${\rm ,} depending
on $\Omega${\rm ,} $A${\rm ,} $p${\rm ,} $q${\rm ,} $\delta$ and $\eps$ {\rm (}\/but not on $\Lambda$
and $\Gamma$\/{\rm )} such that{\rm ,} for every $f\in W^{1,p}(\Omega)$ satisfying
$d^q_A(f)\geq \delta\,
\|f\|_{L^q(\Omega)}$
\begin{equation}\label{t4eq}
\|(\nabla+iA)f\|_{L^p(\Lambda)}+ C\, \|(\nabla+iA)f\|_{L^r(\Gamma)}
\geq \left(\frac1{S^{p,q}_\delta+\eps}+E_A^{p,q}\right)\|f\|_{L^q(\Omega)} \ ,\quad
\end{equation}
where $r=\max\{1,qn/(n+q)\}$ if $1\leq q<\infty${\rm ,} and $S^{p,q}_\delta$ is the optimal constant in
{\rm (\ref{ineqpq}).} For $q=\infty$ {\rm (\ref{t4eq})} holds for any $r>n${\rm ,} and $C$ will also depend on $r$.
\endproclaim

The crucial points to note about (\ref{t4eq}) are the constant 1 in front of
the first term on the left side and the constants $E_A^{p,q}$
and $S^{p,q}_\delta$ on the right side -- which are clearly optimal.
The only unknown constant is $C$. Note that $p>r$ by assumption.

The reader might justly wonder how the volume of $\Gamma$ plays a role
in the error term, as claimed above. The following corollary displays
this dependence; its proof consists just of applying H\"older's
inequality to the second term in (\ref{t4eq}).

\proclaimtitle{Explicit volume dependence}
\specialnumber{1} \proclaim{{C}orollary} \label{C2}
Under the same assumptions as in Theorem~3{\rm ,}
\begin{equation}\label{rem3}
\|(\nabla+iA)f\|_{L^p(\Lambda)}+ C\, |\Gamma|^{1/r-1/p}
\|(\nabla+iA)f\|_{L^p(\Gamma)}
\geq \left(\frac1{S_\delta^{p,q}+\eps}+E_A^{p,q}\right)
\|f\|_{L^q(\Omega)} \ .
\end{equation}
\endproclaim

{\it  Remarks.}
1. A weaker inequality is obtained by substituting $$\|(\nabla+iA)
f\|_{L^p(\Omega)}$$ for $\|(\nabla+iA) f\|_{L^p(\Gamma)}$ in the second
term of (\ref{rem3}). In this way the dependence on
$\Gamma$ is solely through its volume (for any given value of
$\|(\nabla+iA) f\|_{L^p(\Omega)}$).

\vglue4pt 2. 
For $1\leq qn/(q+n) <n$ the exponents of $|\Gamma|$ appearing in
Corollary~\ref{C1} are optimal. This can be seen as follows. If $f$ is
supported in a small ball of volume $|\Gamma|$, the corresponding
minimal \lq energy\rq\ $\|(\nabla+iA) f\|_{L^p}/\|f\|_{L^q}$ is of the
order $|\Gamma|^{1/p-1/q-1/n}$.  Inequality (\ref{rem3}) cannot
hold for a larger exponent since an $f$ supported on several  disjoint
small balls of volume $|\Gamma|$
can be chosen so that $d^q_A(f)\geq
\delta\|f\|_{L^q(\Omega)}$. This would violate (\ref{rem3}) for small
enough $|\Gamma|$. Note that in order to obtain the optimal exponent it is necessary to have Theorem~\ref{T4} in the
critical case $r=qn/(n+q)\geq 1$.

If $q=\infty$ or $q<n/(n-1)$ the optimal dependence on the volume of
$\Gamma$ remains an open problem. In particular this is the case for $n=1$.

\vglue4pt 3. It is  clear that (\ref{t4eq}) cannot hold for $\eps=0$. The
constant $C$ has to go to infinity as $\eps\to 0$. Otherwise, the
inequality would be violated by an $f$ that yields equality in
(\ref{ineqpq}).

\vglue4pt

The proofs of Theorems~\ref{Tmagn} and~\ref{T4} will be given in the
next section. In Section~\ref{a0} we consider the special case $A=0$,
and in Section~\ref{l2sect} we comment on the $L^2$-case $p=q=2$.
Section~\ref{sect5} contains some open problems.

The inequalities in this paper can obviously be extended in various
ways, e.g., to smooth compact manifolds \cite{hebey}, weighted Sobolev
spaces \cite{opic}, or $W^{m,p}(\Omega)$ for $m>1$, but we resist the
temptation to do so here. In fact, in the physics application \cite{lsy02},
Theorem~\ref{T4} is needed for a cube in $\R^3$ with a pair of
opposite faces identified.
 
\vglue-3pt
\section{Proof of Theorems~\ref{Tmagn} and~\ref{T4}}
\advance\eqcount by 8
 \vglue-3pt

{\it Proof of Theorem~2}.
As in the proof of Theorem~\ref{T1} (see~\cite{anal}) we use a
compactness argument. Suppose (\ref{ineqpq}) is false. Then there
exists a sequence of functions $f_j\in W^{1,p}(\Omega)$, with
$\|f_j\|_{L^q(\Omega)}=1$, such that $d^q_A(f_j)\geq \delta$ and
\begin{equation}\label{aaa}
\lim_{j\to\infty} \|(\nabla+iA)f_j\|_{L^p(\Omega)}=E_A^{p,q} \ .
\end{equation}
The sequence $f_j$ is bounded in $L^q(\Omega)$, and it follows from
Theorem~\ref{T1} and (\ref{aaa}) that $f_j$ is actually bounded in
$L^{p}(\Omega)$. Since $A$ is bounded by assumption, $f_j$ is also
bounded in $W^{1,p}(\Omega)$. Hence there exists a subsequence, still
denoted by $f_j$, and a function $f\in W^{1,p}(\Omega)$ such that
$\nabla f_j\rightharpoonup
\nabla f$ weakly in $L^p(\Omega)$. (Note that $p>1$ is important here; for $p=\infty$, weak convergence has to be replaced by weak-* convergence.) The
Rellich-Kondrashov theorem \cite[Thm.~6.2]{adams}
 implies that $f_j\to
f$ strongly in $L^r(\Omega)$ for all $1\leq r<np/(n-p)$ for $p\leq n$,
and for all $1\leq r\leq \infty$ if $p>n$. Hence, by our assumptions
on $q$, $\|f\|_{L^q(\Omega)}=1$ and, by weak lower semicontinuity,
\begin{equation}
\lim_{j\to\infty} \|(\nabla+iA)f_j\|_{L^p(\Omega)}\geq
\|(\nabla+iA)f\|_{L^p(\Omega)}\geq E_A^{p,q}\|f\|_{L^q(\Omega)}= E_A^{p,q}\ . \quad
\end{equation}
This shows that $f\in {\cal  M}_A^{p,q}$ and hence contradicts the fact,
which follows from strong convergence, that
$d_A^q(f)\geq \delta$.
\hfill\qed\vglue8pt 

Before giving the proof of Theorem~\ref{T4}, we state the following
lemma, which is needed to prove Theorem~\ref{T4} in the
critical case $r=qn/(n+q)\geq 1$. It establishes the Poincar\'e inequalities
for functions that vanish on a set of positive measure.

\proclaimtitle{Poincar\'e inequalities for functions with small support}
\specialnumber{1}\proclaim{Lemma}
\label{lem1}
Let $\Omega${\rm ,} $p$ and $q$ be as in Theorem~1{\rm ,} and let
$0<\delta<1$. Then there is a finite number $\widetilde S_{p,q}>0${\rm ,}
which depends on $\Omega,\, \delta,\, p, \, q${\rm ,} such that for any
$f\in W^{1,p}(\Omega)$ with $|\{x\, : \, f(x)\neq 0\}|\leq |\Omega|(1-\delta)$
\begin{equation}\label{delta}
\Vert f \Vert_{L^q(\Omega)} \leq \widetilde S_{p,q} \Vert
\nabla f \Vert_{L^p(\Omega)} \  .  
\end{equation}
\endproclaim

{\it Proof}.
Since $\Omega$ is bounded it suffices to prove this lemma for the
largest possible $q$, given $p$. In particular, it is sufficient to
consider the case $q>1$. From Theorem~\ref{T1} we know that
\begin{equation}\label{14}
\Big\|f-\frac 1{|\Omega|} \int_\Omega f \Big\|_{L^q(\Omega)} 
\leq S_{p,q} \|\nabla f\|_{L^p(\Omega)}
\end{equation}
for the $p$'s and $q$'s in question. By the triangle inequality
\begin{equation}\label{15}
\|f\|_{L^q(\Omega)}\leq \Big\|f-\frac 1{|\Omega|} \int_\Omega f \Big\|_{L^q(\Omega)} + |\Omega|^{1/q-1} \left|\int_\Omega f\right| \ .
\end{equation}
By H\"older's inequality and the assumption on the support of $f$
\begin{equation}\label{16}
\left|\int_\Omega f\right| \leq \|f\|_{L^q(\Omega)} \Big|\{x\, : \, 
f(x)\neq 0\}\Big|^{1-1/q}\leq 
\|f\|_{L^q(\Omega)} \Big[(1-\delta)|\Omega|\Big]^{1-1/q} \ .
\end{equation}
Inserting (\ref{15}) and (\ref{16}) in (\ref{14}) we arrive at
\begin{equation}
\|f\|_{L^q(\Omega)} \leq S_{p,q} \left(1-(1-\delta)^{1-1/q}\right)^{-1} \|\nabla f\|_{L^p(\Omega)} \ .
\end{equation}
\vglue-20pt 
\hfill\qed

\phantom{square}
\vglue-12pt
\demo{Proof of Theorem~3}
Assume that the assertion of the theorem is false. Then there exists a
sequence of triples $(C_j, f_j, \Gamma_j)$, with  $\|f_j\|_{L^q(\Omega)}=1$, $d^q_A(f_j)\break\geq \delta$ and
$\lim_{j\to\infty} C_j=\infty$, such that
\begin{equation}\label{exists}
\lim_{j\to\infty} \left( \|(\nabla+iA)f_j\|_{L^p(\Lambda_j)}
+ C_j  \|(\nabla+iA)f_j\|_{L^r(\Gamma_j)} \right)
< 1/S_\delta^{p,q} +  E_A^{p,q}  \ ,\enspace
\end{equation}
where we set $\Lambda_j=\Omega\setminus \Gamma_j$.
This implies in particular that $\|(\nabla+iA)f_j\|_{L^r(\Gamma_j)}\to 0$ as $j\to\infty$.

We claim that it is no restriction to assume that
$\lim_{j\to\infty} |\Gamma_j|=0$. If this is not the case, define $\gamma_j\subset\Gamma_j$ by
\begin{equation}
\gamma_j=\Big\{x\in \Gamma_j \, : \,  |(\nabla+iA(x))f_j(x)|\geq \|(\nabla+iA)f_j\|_{L^r(\Gamma_j)}^{1/2} \Big\} \ .
\end{equation}
Note that $|\gamma_j|\leq  \|(\nabla+iA)f_j\|_{L^r(\Gamma_j)}^{r/2}\to 0$ as $j\to \infty$. Moreover,
\begin{equation}
\|(\nabla+iA)f_j\|_{L^p(\Gamma_j\setminus \gamma_j)} \leq |\Omega|^{1/p}
\|(\nabla+iA)f_j\|_{L^r(\Gamma_j)}^{1/2} \ ,
\end{equation}
which also goes to zero as $j\to\infty$. Therefore (\ref{exists})
holds with $\Lambda_j$ and $\Gamma_j$ replaced by $\Lambda_j \cup
(\Gamma_j\setminus\gamma_j)$ and $\gamma_j$, respectively.

It suffices, therefore, to consider the case $\lim_{j\to\infty}
|\Gamma_j|=0$. By passing to a
subsequence we can assume that $\sum_j |\Gamma_j|$ is finite.

The sequence $f_j$ is bounded in $L^q(\Omega)$, and also
$(\nabla+iA)f_j$ is bounded in $L^r(\Omega)$. From Theorem~\ref{T1} we
see that $f_j$ is actually bounded in $L^{\tilde{q}}(\Omega)$ for
$\tilde{q}=\max\{q,n/(n-1)\}>1$, and hence there is a subsequence,
still denoted by $f_j$, and an $f\in L^q(\Omega)$, such that
$f_j\rightharpoonup f$ weakly in $L^{q}(\Omega)$. (If $q=\infty$ weak
convergence has to be replaced by weak-* convergence.)

For some fixed $N$ let $\Sigma_N=\Omega\setminus \bigcup_{j\geq N}
\Gamma_j$. Note that $(\nabla+iA)f_j$ is bounded in
$L^p(\Sigma_N)$ and, therefore, we can choose a subsequence such
that $(\nabla+iA)f_j\rightharpoonup (\nabla+iA)f$ weakly in
$L^p(\Sigma_N)$. (Again, replace weak by weak-* if  $p=\infty$.) By weak lower semicontinuity of norms and the
fact that $\Sigma_N\subset \Lambda_j$ for $j\geq N$,
\begin{equation}\label{semic}
\liminf_{j\to\infty} \|(\nabla+iA)f_j\|_{L^p(\Lambda_j)}
\geq \liminf_{j\to\infty}  \|(\nabla+iA)f_j\|_{L^p(\Sigma_N)}
\geq \|(\nabla+iA)f\|_{L^p(\Sigma_N)} \ .
\end{equation}
This holds for all $N$ and, since $\Sigma_{N}\subset
\Sigma_{N+1}$ and $|\bigcup_N\Sigma_N|=|\Omega|$,
\begin{equation}\label{semic2}
\liminf_{j\to\infty} \|(\nabla+iA)f_j\|_{L^p(\Lambda_j)}
\geq \|(\nabla+iA)f\|_{L^p(\Omega)} \ .
\end{equation}
Suppose that we knew that $f_j\to f$ {\it strongly} in
$L^q(\Omega)$.  Then clearly $d_A^q(f)\geq \delta
\|f\|_{L^q(\Omega)}=\delta$, so the right side of (\ref{semic2})
would be $\geq 1/S_\delta^{p,q}+ E_A^{p,q}$ by (\ref{ineqpq}), and
thereby contradict (\ref{exists}) and establish (\ref{t4eq}).

In the following, we will show that $f_j\to f$ {\it strongly} in
$L^{q}(\Omega)$. Note that for $q<n/(n-1)$ and for $q=\infty$ this follows
immediately from the Rellich-Kondrashov theorem
\cite[Thm.~6.2]{adams},
so we can restrict ourselves to the case $1<{q}<\infty$. For $M>0$, define
\begin{equation}
f_j^M(x) = \min\{M, |f_j(x)|\}
\end{equation}
and
\begin{equation}
h_j^M(x)=|f_j(x)|-f_j^M(x) \ .
\end{equation}
Note that both $f_j^M$ and $h_j^M$ are in $W^{1,p}(\Omega)$. Moreover,
\begin{equation}\label{ref3}
\Big|\Big\{x\, :\, h_j^M(x)\neq 0\Big\}\Big|=\Big |\Big\{x\in\Omega\, : \, |f_j(x)|>M\Big\}\Big| \leq \frac
{\|f_j\|_{L^q(\Omega)}^q}{M^q}=\frac 1{M^q} \ .\quad
\end{equation}
By choosing $M$ larger that $2|\Omega|^{-1/q}$ we can use
Lemma~\ref{lem1} to conclude that
\begin{equation}\label{ref1}
\|h_j^M\|_{L^{q}(\Omega)}\leq S\|\nabla h_j^M\|_{L^{r}(\Omega)}
\end{equation}
for some constant $S$ independent of $M$ and $j$. Note that the
intersection of the two sets
$\alpha_j:=\{x\, :\, \nabla
f_j^M(x)\neq 0\}$ and $\beta_j:=\{x\, :\, \nabla h_j^M(x)\neq 0\}$ has
measure zero.  Therefore
\begin{equation} \label{ref2}
\|\nabla h_j^M\|_{L^{r}(\Omega)}=\|\nabla |f_j|\|_{L^{r}(\beta_j)}
\leq \|(\nabla+iA)f_j\|_{L^p(\Lambda_j)}|\beta_j|^{1/r-1/p} +
\|(\nabla+iA)f_j\|_{L^r(\Gamma_j)} \ ,
\end{equation}
where we used again H\"older's inequality and also the diamagnetic
inequality $|\nabla |f|(x)|\leq |(\nabla+iA(x))f(x)|$ (see
\cite[Thm.~7.21]{anal}). By (\ref{ref3}), $|\beta_j| \leq 1/M^q$.
This fact, together with (\ref{ref1}),
(\ref{ref2}),  and (\ref{exists}), implies that
\begin{equation}\label{hres}
\limsup_{j\to\infty} \|h_j^M\|_{L^{ q}(\Omega)}\leq S (E_A^{p,q}+1/S_\delta^{p,q}) M^{(r-p)/prq} \ .
\end{equation}
(If $p=\infty$, the exponent in the last term has to be replaced by $-1/rq$.)

{}From (\ref{exists}) we see that $(\nabla+iA)f_j$ is a bounded
sequence in $L^{r}(\Omega)$ and, since $A$ is bounded by assumption,
the same is true for $\nabla f_j$. Hence we can apply the
Rellich-Kondrashov theorem (see, e.g., \cite[Thm.~6.2]{adams}) to conclude that, modulo choice of a
subsequence, $f_j\to f$ {\it strongly} in $L^{ q-\nu}(\Omega)$ for any
$0<\nu\leq  q-1$, and therefore
\begin{equation}
\int_\Omega |f|^{ q-\nu} = \lim_{j\to\infty} \int_\Omega |f_j|^{ q-\nu} \ .
\end{equation}
By definition of $f_j^M$,
\begin{equation}
\int_\Omega |f_j|^{ q-\nu} \geq \int_\Omega |f_j^M|^{ q-\nu}\geq \frac 1{M^\nu} \int_\Omega |f_j^M|^{ q} \ .
\end{equation}
Using (\ref{hres}) we therefore obtain
\begin{equation}
\int_\Omega |f|^{ q-\nu}\geq  \frac 1{M^\nu} \left(1 - \big[S (E_A^{p,q}+1/S_\delta^{p,q})\big]^{ q} M^{ (r-p)/pr}\right) \ ,
\end{equation}
and hence
\begin{equation}
\int_\Omega |f|^{ q} = \lim_{\nu\to 0} \int_\Omega |f|^{ q-\nu} \geq 1 - \big[S (E_A^{p,q}+1/S_\delta^{p,q})\big]^{ q} M^{ (r-p)/pr}  \ .
\end{equation}
Since $M$ can be chosen arbitrarily large, and $p>r$, this shows that
\hbox{$\|f\|_{L^{ q}(\Omega)}\!=\!1$,} implying strong convergence and finishing
the proof.
\enddemo

As might be expected, the proof of (\ref{t4eq}) can be simplified if
one is not interested in the optimal $r$, but rather
$r>\max\{1,qn/(n+q)\}$.

\section{The special case 
$A=0$}
\label{a0}
\advance\eqcount by 30

In the case of vanishing magnetic field $A=0$, there is a much simpler
proof of Theorem~\ref{T4}. In fact this theorem follows easily from
Theorem~\ref{T1}, as we now show. However, this simple proof has the
disadvantage of not yielding any information about the optimal
constants.

\proclaimtitle{Generalized  Poincar\'e inequalities for
$A=0$} \specialnumber{4} \proclaim{Theorem} \label{T2} Let $\Omega${\rm ,} $g${\rm ,} $p${\rm ,} $q$ be as in
Theorem~1{\rm ,} and let $\tilde{q}_n=\max\{1,qn/(n+q)\}$. Let
$\Lambda
\subset \Omega$ be a measurable subset of $\Omega$ and let
$\Gamma=\Omega\setminus\Lambda$.

There are constants $S^{p,q}$ {\rm (}\/generally different from $S_{p,q}${\rm ),}
depending only on $\Omega${\rm ,} $g${\rm ,} $p$ and $q${\rm ,} but not on $\Lambda${\rm ,}
such that  for all $f \in W^{1,p}(\Omega)$
\begin{equation}  \label{firstineq}
\Big\Vert f- \int_{\Omega}fg\, \Big\Vert_{L^q(\Omega)} \leq S^{p,q}\Big[ \Vert
\nabla f \Vert_{L^p(\Lambda)} + |\Omega|^{1/p-1/\tilde{q}_n}\Vert
\nabla f \Vert_{L^{\tilde{q}_n}(\Gamma)}  \Big]  \hskip.5in
\end{equation}
if $1\leq q< \infty$ and $\tilde{q}_n\leq p\leq \infty$.
 One can take $S^{p,q}=S_{\tilde{q}_n,q}|\Omega|^{1/\tilde{q}_n-1/p}$.

For $q=\infty${\rm ,} there exist  constants $\widehat S^{p,r}$ such that
\begin{equation}  \label{firstineqpprime}
\Big\Vert f- \int_{\Omega}fg\, \Big\Vert_{L^\infty(\Omega)} \leq \widehat
S^{p,r}\Big[ \Vert
\nabla f \Vert_{L^p(\Lambda)} + |\Omega|^{ 1/p-1/r}\Vert
\nabla f \Vert_{L^r(\Gamma)} \Big]
\end{equation}
for all $n<r\leq p\leq \infty$. Now{\rm ,} $\widehat
S^{p,r}=S_{p,\infty}|\Omega|^{1/r-1/p}$.
\endproclaim

As in Corollary~\ref{C2}, the application of H\"older's inequality to
the rightmost norms in (\ref{firstineq}) and (\ref{firstineqpprime})
displays the dependence on the volume of $|\Gamma|$. We obtain

 \proclaimtitle{Explicit volume dependence}
\specialnumber{2} \proclaim{{C}orollary}\label{C1}
Under the assumptions of Theorem~4
\begin{equation}  \label{firstineqe}
\Big\Vert f- \int_{\Omega}fg\, \Big\Vert_{L^q(\Omega)} \leq S^{p,q}\left[ \Vert
\nabla f \Vert_{L^p(\Lambda)} + \left(\frac{|\Gamma|}{|\Omega|}
\right)^{1/\tilde{q}_n - 1/p}\Vert
\nabla f \Vert_{L^p(\Gamma)} \right]  \hskip.25in
\end{equation}
if $1\leq q< \infty$ and $\tilde{q}_n\leq p\leq \infty$.

For $q=\infty${\rm ,}
\begin{equation}  \label{firstineqpprimee}
\Big\Vert f- \int_{\Omega}fg\, \Big\Vert_{L^\infty(\Omega)} \leq \widehat
S^{p,r}\left[ \Vert
\nabla f \Vert_{L^p(\Lambda)} + \left(\frac{|\Gamma|}{|\Omega|}
\right)^{ 1/r- 1/p}\Vert
\nabla f \Vert_{L^p(\Gamma)} \right]  \hskip.25in
\end{equation}
for all $n<r\leq p\leq \infty$.
\endproclaim

{\it Remarks.}
1. As a special case, we can assume that $g(x)=0$ for $x\in \Gamma$ in
Corollary~\ref{C1}, and use the simple fact that $\Vert \cdot
\Vert_{L^q(\Omega)} \geq
\Vert \cdot \Vert_{L^q(\Lambda)}$ to obtain
\begin{equation}  \label{firstineqprime2}
\Big\Vert f- \int_{\Lambda}fg\, \Big\Vert_{L^q(\Lambda)} \leq S^{p,q} \left[ \Vert
\nabla f \Vert_{L^p(\Lambda)} +\left(\frac{|\Gamma|}{|\Omega|}
\right)^{1/{\tilde{q}_n} - 1/p}\Vert
\nabla f \Vert_{L^p(\Gamma)}  \right]  \hskip.25in
\end{equation}
when $q<\infty$, and similarly for (\ref{firstineqpprime}). The virtue of (\ref{firstineqprime2}) is that it is an inequality that depends only on $\Lambda$, except for an error term.
We emphasize again  that the constants $S^{p,q}$ do
not depend on $\Lambda$, but only on $\Omega$ and $g$.

\vglue8pt 2.
Theorem~\ref{T2} is a corollary of Theorem~\ref{T1}. This is in contrast
to our general result, Theorem~\ref{T4}, which does not
appear to follow easily from Theorem~\ref{T1}.

\vglue8pt 3.
The optimal constant $S^{p,q}$ in Theorem~\ref{T2} is left
unspecified. This is in contrast to Theorem~\ref{T4}, where the
constant appearing on the right side of (\ref{t4eq}) is optimal, up to
an~$\eps$. The simple proof we shall give of Theorem~\ref{T2}, as a
corollary of Theorem \ref{T1}, does not
allow us to relate $S^{p,q}$ to the optimal constant for the usual
Poincar{\'e} inequality for $\Lambda=\Omega$ (although we can relate
it to $S_{\tilde q_n,q}$). Thus, even in the $A=0$ case, the more
complicated proof of Theorem~\ref{T4} has the advantage of yielding
information about the sharp constant.

\vglue12pt 4. 
In contrast to Theorem~\ref{T4}, Theorem~\ref{T2} includes the
critical case $p=\tilde{q}_n$. Note, however, that in this case
Theorem~\ref{T2} does not represent any improvement over
Theorem~\ref{T1}.

\demo{Proof of Theorem~4}
By Theorem~\ref{T1} with $1\leq q <\infty$,
\begin{equation} \label{poinn}
\Big\Vert f- \int_{\Omega}fg \Big\Vert_{L^q(\Omega)} \leq
S_{\tilde{q}_n,q} \Vert
\nabla f \Vert_{L^{\tilde{q}_n}(\Omega)} \  .
\end{equation}
We estimate the right side by the triangle inequality
\begin{equation}
\|\nabla f\|_{L^{\tilde{q}_n}(\Omega)}\leq \|\nabla f\|_{L^{\tilde{q}_n}(\Lambda)} +
\|\nabla f\|_{L^{\tilde{q}_n}(\Gamma)} \ .
\end{equation}
H\"older's inequality implies that, for any $p\geq \tilde{q}_n$,
\begin{equation}
\|\nabla f\|_{L^{\tilde{q}_n}(\Lambda)}\leq \|\nabla f\|_{L^p(\Lambda)} |\Lambda|^{1/\tilde{q}_n-1/p} \ .
\end{equation}
By the fact that  $|\Lambda|\leq
|\Omega|$, this proves (\ref{firstineq}), with
$S^{p,q}=S_{\tilde{q}_n,q}|\Omega|^{1/\tilde{q}_n-1/p}$.
The same proof works for $q=\infty$, with $\widehat S^{p,r}=S_{p,\infty}|\Omega|^{1/r-1/p}$.
\enddemo

\section{The special case   
$p=q=2$}
\label{l2sect}
\advance\eqcount by 38

In the case $p=q=2$, the ground state manifold ${\cal  M}_A^{2,2}$
is a linear subspace of $L^2(\Omega)$, spanned by the eigenfunctions
corresponding to the lowest (Neumann) eigenvalue of the operator
$H=-(\nabla+iA)^2$. And Theorem~\ref{Tmagn} is just the statement that
there is a gap above the lowest eigenvalue, which follows from the
discreteness of the spectrum of $H$. The dimension of ${\cal 
M}_A^{2,2}$ is finite, but it can be strictly greater than one. These properties are shown in the appendix. This is in contrast to the case $A=0$, where
${\cal  M}_A^{2,2}$ is one-dimensional, spanned by the constant
function.

We can use Theorem~\ref{T4} in the case $p=q=2$ to get an inequality
of the form (\ref{l5}) that holds for all $f\in W^{1,2}(\Omega)$. For
simplicity we state it for the case of ${\cal  M}_A^{2,2}$ being
one-dimensional. The proof
is obtained by replacing $f$ in (\ref{rem3}) by $f-\phi_A \int_\Omega
f g$ and using the Cauchy-Schwarz inequality.

\proclaimtitle{analogue of (\ref{l5}) for punctured domains}
\specialnumber{3} \proclaim{{C}orollary} \label{C3}
Let
$\Omega \subset \R^n$ be a bounded{\rm ,} connected{\rm ,} open set that has the
cone property{\rm ,} and let $E_A^{2,2}$ and ${\cal  M}_A^{2,2}$ be as explained
above. Let $\Lambda
\subset \Omega$ be a measurable subset of $\Omega${\rm ,}
 $\Gamma=\Omega\setminus\Lambda${\rm ,} and let $0<\delta\leq 1$.
Suppose that ${\cal  M}_A^{2,2}$ is one\/{\rm -}\/dimensional{\rm ,} spanned by the
normalized eigenfunction $\phi_A$ corresponding to $E_A^{2,2}$. Let
$g$ be an $L^2(\Omega)$ function \pagebreak satisfying $\int_{\Omega} g\phi_A
=1$. For any $\eps>0$ there exists a positive constant $C'$ depending
on $\Omega${\rm ,} $A${\rm ,} $g$ and $\eps$ {\rm (}\/but not on $\Lambda$ and $\Gamma${\rm )}
such that{\rm ,} for every $f\in W^{1,2}(\Omega)${\rm ,}
\begin{eqnarray}  \label{cor}
&&\hskip-.25in \|(\nabla+iA) f\|_{L^2(\Lambda)}^2+C'\,
|\Gamma|^{\min\{1,2/n\}}\|(\nabla+iA)f\|_{L^2(\Gamma)}^2 \\ &&\hskip.95in +  C'\, \Big| \int_\Omega
fg \Big|\, \|(\nabla+iA)f\|_{L^2(\Gamma)}
\|(\nabla+iA)\phi_A\|_{L^2(\Gamma)}\nonumber \\
&&\hskip.8in\geq \big(E_A^{2,2}\big)^2 \|f\|_{L^2(\Omega)}^2 + \frac 1{S^g+\eps}
\Big\|f -\phi_A \mbox{$\left[\int_\Omega fg \right]$}
\Big\|_{L^2(\Omega)}^2\ ,\nonumber
\end{eqnarray}
with $S^g$   related to the optimal constant $\widetilde S_g^{2,2}$ in {\rm (\ref{l5})} by
$$S^g=  \big(S_g^{2,2}\big)^2\big(1+ 2 E_A^{2,2} S_g^{2,2}\big)^{-1}.$$
\endproclaim

\vglue-12pt
{\it Remarks.}
1. By the same argument as in Remark~2 after Theorem~\ref{T4} the
exponent $2/n$ in (\ref{cor}) is sharp for $n\geq 2$; i.e., it cannot
be increased. It is natural to conjecture that (\ref{cor}) holds with
exponent $2=2/n$ also for $n=1$.  Unfortunately, the method of proof
presented here does not allow for this generalization.

\vglue4pt 2.
For regular enough boundary of $\Omega$ it follows from
elliptic regularity that $(\nabla+iA)\phi_A$ is in fact a
bounded function \cite{gilbarg}. This allows us to replace
$ \|(\nabla+iA)\phi_A\|_{L^2(\Gamma)}$
by $\const |\Gamma|^{1/2} $. In any
case, $ \|(\nabla+iA)\phi_A\|_{L^2(\Gamma)}$
goes to zero as $|\Gamma|\to 0$.

\vglue4pt 3.
As in Remark 1 after Theorem~\ref{T2} one can consider the special
case where $g$ vanishes on $\Gamma$ to obtain an inequality that depends on $\Gamma$ only via its volume. It has to be noted, however, that $E_A^{2,2}$ is defined on the whole of $\Omega$ and not just on $\Lambda$.

\vglue4pt 4.
Strictly speaking, Corollary~\ref{C3} is not really a corollary of
Theorem~\ref{T4} because of the optimal constant $S^g$ appearing in
(\ref{cor}). From (\ref{rem3}) we can only infer (\ref{cor}) with $S^g$
replaced by $S_\delta$, for some $\delta$ depending on $g$.  However,
by imitating the proof of Theorem~\ref{T4} one can show that (\ref{cor})
holds.
 
\vglue-5pt
\section{Some open problems}\label{sect5}
\vglue-4pt

$\bullet$\enspace In Theorems~\ref{Tmagn} and~\ref{T4} we have excluded the `critical'
case $p=\max\{1,qn/\break (n+q)\}$. In this case, the existence of minimizers
of (\ref{defea2}) and hence the nonemptiness of ${\cal 
M}_A^{p,q}$ is {\it a priori} not clear, except for the case $A=0$,
where ${\cal  M}_A^{p,q}$ is trivially just the one-dimensional
space spanned by the constant function.

\vglue4pt $\bullet$\enspace 
The optimal exponent in the dependence on the volume of $\Gamma$ in
(\ref{rem3}) (see Remark~2 after Corollary~\ref{C2}) is open for the
cases $q=\infty$ and $q<n/(n-1)$. This comes from the fact that
(\ref{rem3}) is obtained as a corollary of (\ref{t4eq}), where
necessarily $r\geq 1$.

  $\bullet$\enspace 
Since Theorem~\ref{T4} is proved by a compactness argument, the
constant $C$ appearing in (\ref{t4eq}) is left unspecified. It would
be desirable to obtain a decent upper bound for this value. Also the
optimal values of the constants $S_\delta^{p,q}$ appearing in
Theorem~\ref{Tmagn} are in general unknown. Indeed, decent estimates
of $S_{p,q}$ in Theorem~\ref{T1} are not readily available when $p,q \neq 2$.

\vglue4pt $\bullet$\enspace 
The dimension of the ground state manifold ${\cal  M}_A^{p,q}$ can
be bigger that one. In the appendix we give an example for the case
$p=q=2$ where its dimension is two. It would be interesting to
construct examples (or prove their existence) where the dimension can
be arbitrarily large.

\advance\eqcount by 38

\vglue22pt \centerline {\bf Appendix: Spectrum of 
$-(\nabla +iA)^2$}
\vglue16pt 
Here we prove two facts about $-(\nabla +iA)^2$  which were used in the
text. As before, we have a bounded, connected domain $\Omega$ in $\R^n$
with the cone property. Connectedness is not really necessary here,
but the number of connected components should be finite. The vector
field $A$ is bounded and measurable. The boundedness is not crucial
but we assume it for simplicity.

We define the eigenvalues (spectrum) $E_k$ and eigenfunctions $\phi_k$
of\break $-(\nabla +iA)^2$ in $L^2(\Omega)$ by means of quadratic forms as
in \cite{anal}, i.e.,  $E_{k+1}$ is defined by
\begin{equation}\nonumber
E_{k+1} =  \inf \left\{\Vert
(\nabla +i A)f\Vert_{L^2(\Omega)}^2  \, : \, \Vert f\Vert_{L^2(\Omega)}^2
=1 , \, \int_\Omega f \overline {\phi_j} =0 \ {\rm for}\ j=1,\ldots ,k \right\}  \ .
\end{equation}
Then, by standard methods (using the Rellich-Kondrashov theorem), one
shows that there is a minimizer for $E_{k+1}$, which is called
$\phi_{k+1}$.  In the text, $\sqrt{E_1}$ was called $E_A^{2,2}$.

The two facts are the following.

\proclaimtitle{Spectrum of   
$-(\nabla +iA)^2$} 
\specialnumber{2}\proclaim{Lemma}\label{spec} 
\begin{itemize}
\item[{\rm A.}]
The spectrum is discrete\/{\rm ;} i.e.{\rm ,} the number of
eigenvalues less than any number $E$ is finite.
\item[{\rm B.}]
The multiplicity of $E_1$ can be greater than one.
\end{itemize}

\endproclaim

\demo{Proof}
To prove A we suppose that there are infinitely many eigenvalues below
$E$. If $\psi \in W^{1,2}(\Omega)$, with $\|\psi\|_{L^2(\Omega)}=1$ and $\Vert (\nabla
+iA)\psi\Vert_{L^2(\Omega)} \leq E^{1/2}$, then, since $A$ is bounded,
$\Vert \nabla \psi\Vert_{L^2(\Omega)}\leq E^{1/2}+\|A\|_{L^\infty(\Omega)}$.  Thus, the infinite
sequence of functions $\phi_1,\, \phi_2, ...$ is bounded in
$W^{1,2}(\Omega)$. By the Rellich-Kondrashov theorem this sequence has
a subsequence that converges strongly in $L^2(\Omega)$. This is
impossible since the $\phi_i$'s are orthonormal (and hence\break $\Vert
\phi_i-\phi_j\Vert_{L^2(\Omega)} = \sqrt 2$ for $i\neq j$).

To prove B, consider the case in which $\Omega\subset \R^2$ is an
annulus centered at the origin (or a cylinder in $\R^n$ based on an
annulus in $\R^2$). We take $A$ to be (in polar coordinates)
$A(r,\theta)=\lambda r^{-1} \widehat e_\theta$, where $\widehat
e_\theta $ is the unit vector in the $\theta$ direction. We shall
show that for suitable values of $\lambda$ the multiplicity of $E_1$
is two.

The Hilbert space $W^{1,2}(\Omega)$ is the direct sum of Hilbert
spaces $W^{1,2}_l(\Omega)$, $l\in {\Bbb Z}$, consisting of
functions of the form $\exp(-il\theta)g(r)$, and these subspaces
continue to be orthogonal under the action of $\nabla +iA$. Thus, the
eigenvectors in our case can be chosen to belong to exactly one of
these subspaces. Thus, we can define $E_A(l)$ to be the lowest
eigenvalue in $W^{1,2}_l(\Omega)$, and $E_1$ is then the minimum among
the numbers $E_A(l)$.  Since
\begin{equation}
-(\nabla +iA)^2=-\frac {\partial^2}{\partial r^2}-\frac 1r \frac
{\partial}{\partial r}-\frac 1{r^2}\left( \frac {\partial}{\partial
\theta} + i\lambda\right)^2 \ ,
\end{equation}
$E_1=E_A(l_0)$, where $l_0$ is the integer closest to
$\lambda$. Therefore, for $\lambda \in {\Bbb Z}+\half$, there are
two eigenfunctions with the same eigenvalue
$E_1=E_A(\lambda-\half)=E_A(\lambda+\half)$. \phantom{second}
\enddemo

\end{document}